\documentclass[article,12pt]{elsarticle}

\makeatletter
\def\ps@pprintTitle{%
 \let\@oddhead\@empty
 \let\@evenhead\@empty
 \def\@oddfoot{\centerline{\thepage}}%
 \let\@evenfoot\@oddfoot}
\makeatother
\usepackage{amssymb}

\journal{Indagationes Mathematicae}

\begin{document}

\begin{frontmatter}

\newtheorem{theorem}{Theorem} 
\newtheorem{lemma}[theorem]{Lemma}
%%\newtheorem{theorem}{Lemma}
%%\theoremstyle{definition}
%%\newtheorem{defn}[thm]{Definition}



\title{Modular Arithmetic Study on Ascending and Descending Operations Related with Collatz Conjecture and Origin of the Operation}


\author{Kyo Jin Ihn}

\affiliation{organization={Department of Chemical Engineering, Kangwon National University},%Department and Organization
 %%           addressline={}, 
            city={Chunchon},
            postcode={24341}, 
 %%           state={},
            country={Korea (South)}}

\begin{abstract}
Sequence of numbers generated by the recurrence relation based on the Collatz conjecture is investigated. An arithmetic operation on the Collatz conjecture is called descending operation, and ascending operation is carried out reversely to the descending operation. Study on the $3x+1$ problem against every integer is relevant to that of descending operation against every odd. Any odds can be generated by ascending operations against odds not multiple of 3, and vise versa. Ascending operations against an odd not multiple of 3 can generate infinite number of odds in the next generation. The odds multiple of 3 are terminal numbers of sequences that no further ascending operations are possible. Descending operations against the odds 1,5,21, $\cdots \left(= \sum_{j=0}^{m} 4^j \right)$ reach 1, and otherwise reach the odds greater than or equal to 5. No two odds except 1 in each sequence are the same, and the each descending sequence is unique. Every descending sequence can start from terminal numbers, and reaches ultimately the origin of 1 as indicated by the Collatz conjecture. Sequences of odds in the same pattern ignoring the sizes of the odds are generated by ascending operations in the same number of doubling operations against the odds in modular arithmetic.
\end{abstract}


\end{frontmatter}

%%===================================================================
\bigskip
%%===================================================================


\section{Introduction}
\label{}


The Collatz conjecture introduced in 1937 by Lothar Collatz examines whether repeating arithmetic operations will transform every positive integer into 1.\cite{CollatzBook} It is known as 3x+1 problem.\cite{CollatzBook} The 3x+1 problem involves the operations generating a sequence of numbers as follows: $3x+1$ operation if it is odd, and halving operation if it is even. The conjecture is that the sequences generated by the consecutive operations against any positive integers reach always 1.\cite{CollatzBook} 

By the Collatz conjecture, halving operations until to reach an odd followed by $3x+1$ operation against a positive odd $x$ is called descending operation in this study. If we let $x_{k-1}$ be a positive odd generated by a descending operation against a positive odd $x_k$, number of halving operations is determined by p-adic valuation such that $m_{k-1}=\nu_2 \left( 3x_k +1 \right)$, .

The descending operation includes two operations such that $3x_k+1$ operation against an odd $x_k$ and halving operations against the even $3x_k +1$ until to reach an odd $x_{k-1}$. If an even generates an odd $x_k$ by halving operations, it belongs to the evens generated in the process of descending operation against an odd $x_{k+1}$. Every odd is generated by the halving operation against every even. Consequently, in proving the Collatz conjecture, the study on the 3x+1 problem against every positive integer is relevant to the study using the recurrence relation of the descending operation against every positive odd. 

Consecutive descending operations using the recurrence relation $x_{k-1}=\left(3x_k +1\right) / 2^{m_{k-1}}$ generate a sequence of odds. N-terms of descending sequence can be expressed as $x_n\rightarrow x_{n-1}\rightarrow \cdots \rightarrow x_k \rightarrow \cdots \rightarrow x_1\rightarrow x_0$. The $x_0$ is not always the origin of 1 unless it is defined. The odds $x_{n-1}$ to $x_0$ in the sequence are not multiples of 3, since the odds multiples of 3 cannot be generated by the descending operation.

Each sequence of odds generated by consecutive descending operations is unique, because numbers of halving operations are fixed. Descending sequences can be thus expressed by the odds and the numbers of halving operations instead of showing whole evens for simplicity. Collatz conjecture states that every descending sequence reaches ultimately the origin of 1.\cite{CollatzBook} 

The reverse operation against the descending operation such that $\left( 2^m x-1\right)/3$ operation against a positive odd $x$ not multiple of 3 is called ascending operation, where $m$ is number of doubling operations to show $2^m x\equiv 1 \left( mod \: 3 \right)$. The odds generated by the ascending operations from an odd correspond to the first generation odds. Ascending operations against an odd not multiple of 3 can generate infinite number of odds in the first generation, since multiple doubling operations are possible.

\section{Types of ascending operations}
\label{}
\subsection{First generation odds}

Odds $x_k$ are generated by the recurrence relation $x_k= \left( 2^{m_{k-1}}  x_{k-1}-1\right)/3$ of ascending operation with $m_{k-1}$ times of doubling operations against a positive odd $x_{k-1}$ not multiple of 3. Doubling operations against every odd can generate every even. Since $2^2  \equiv 1  \left( mod \: 3 \right)$, $2^{2m} \equiv 1 \left( mod \: 3 \right)$ and $2^{2m-1} \equiv -1 \left( mod \: 3 \right)$, where $m$ is positive integer. $2^{2m}x_0 \equiv 1 \left( mod \: 3 \right)$ and $2^{2m-1}x_0 \equiv -1 \left( mod \: 3 \right)$ by the odds $x_0 \equiv 1 \left( mod \: 3 \right)$, and  $2^{2m-1}x_0 \equiv 1 \left( mod \: 3 \right)$ and $2^{2m}x_0 \equiv -1 \left( mod \: 3 \right)$ by $x_0 \equiv -1 \left( mod \: 3 \right)$. Thus, ascending operations are performed by every two doubling operations. This proves Theorem 1. 

\begin {theorem}
$2^{2m}x \equiv 1 \left( mod \: 3 \right)$ by positive odd $x \equiv 1 \left( mod \: 3 \right)$ and $2^{2m-1}x \equiv 1 \left( mod \: 3 \right)$ by $x \equiv -1 \left( mod \: 3 \right)$, where m is positive integer.
\end {theorem} 

Tables 1 and 2 are the examples of  first generation odds generated by ascending operations concerning up to 12 times of doubling operations against $x_0 = 1$ and 5. Ascending operations against $x_0 =1$ using the recurrence relation $x_1 = \left( 2^m  x_0-1\right)/3$ generate the odds $x_1=\sum_{j=0}^{(m-2)/2} 4^j$ for even $m$  by Theorem 1. Ascending operations against $x_0=5$ generate the odds $x_1=2^m+\sum_{j=0}^{(m-1)/2} 4^j$  for odd $m$.

%% Table 1
\begin{table} [p!]
\centering
\caption {Examples of first generation odds $x_1=\sum_{j=0}^{(m-2)/2} 4^j$ for even $m$ generated by the recurrence relation $x_1 = \left( 2^m  x_0-1\right)/3$ of ascending operation against $x_0 =1$}
\label{table:4}
\centering
\begin{tabular}{c c c c c c} \\ [0.5ex]
\hline
$x_0$ & $m$ & $2^m x_0$ & $x_1$ & $x_1\; (mod\; 3)$ & $2^mx_0\; (mod\; 3)$    \\
  \hline 
1&1 & $2$ & & & \\ 
1& 2 & $2^2$ & 1 & 1 & 1 \\ 
1& 3 & $2^3$ & & & \\ 
1& 4 & $2^4$ & 5 & -1 & 1 \\ 
1& 5 & $2^5$ & & & \\ 
1& 6 & $2^6$ & 21 & 0 & 1 \\ 
1& 7 & $2^7$ & & & \\ 
1& 8 & $2^8$ & 85 & 1 & 1 \\ 
1& 9 & $2^9$ & & & \\ 
1& 10 & $2^{10}$ & 341 & -1 & 1 \\
1& 11 & $2^{11}$ & & & \\ 
1& 12 & $2^{12}$ & 1365 & 0 & 1 \\ 
 \hline
\end{tabular}
\end{table}

%% Table 2
\begin{table} [p!]
\centering
\caption {Examples of first generation odds $x_1=2^m+\sum_{j=0}^{(m-1)/2} 4^j$ for odd $m$ generated by the recurrence relation $x_1 = \left( 2^m  x_0-1\right)/3$ of ascending operation against $x_0 =5$}
\label{table:4}
\centering
\begin{tabular}{c c c c c c} \\ [0.5ex]
\hline
$x_0$ & $m$ & $2^m x_0$ & $x_1$ & $x_1\; (mod\; 3)$ & $2^mx_0\; (mod\; 3)$    \\
  \hline 
5& 1 & $5\times2$ & 3 & 0 & 1 \\ 
5& 2 & $5\times2^2$ & & & \\ 
5& 3 & $5\times2^3 $ & 13 &1 &1  \\ 
5& 4 & $5\times2^4$ & & & \\ 
5& 5 & $5\times2^5$ &53 &-1 &1 \\ 
5& 6 & $5\times2^6$ & & & \\ 
5& 7 & $5\times2^7$ &213 &0 &1 \\ 
5& 8 & $5\times2^8$ & & & \\ 
5& 9 & $5\times2^9$ &853 &1 &1 \\ 
5& 10 & $5\times2^{10}$ & &  & \\
5& 11 & $5\times2^{11}$ &3413 &-1 &1 \\ 
5& 12 & $5\times2^{12}$ & & & \\ 
 \hline
\end{tabular}
\end{table}

If we assume that ascending operations against two odds $x_{k-1}$ and $x_{k-1}+a$ with even $a$ generate the same odd $x_k$, we can let the recurrence relations as $x_k = \left(2^n x_{k-1}-1\right)/3$ and $x_k = \left(2^m \left(x_{k-1}+a\right)-1\right)/3$ with $n>m$. The resulting formula $\left(2^{n-m}-1\right) x_{k-1}=a$ is not correct, since left side is an odd. This indicates that no two odds $x_k$ are the same, if the odds in the previous generation are different. This proves Lemma 2. 

\begin {lemma}
 Let $x_k$ be the odds generated by $x_k = \left(2^{m_{k-1}} x_{k-1}-1\right)/3$ operation against odds $x_{k-1}$ not multiple of 3, then no two odds $x_k$ in the $k$-th generateion are the same, if the odds $x_{k-1}$ in the previous generation are different. 
\end {lemma}

Every odd not multiple of 3 is defined by $x_{k-1}= 6n+1$ and $x_{k-1}= 6n+5$, where $n$ is nonnegative integer. In the recurrence relation $x_k=\left(2^{m_{k-1}} x_{k-1} -1\right)/3$ of ascending operation, if $x_{k-1}= 6n+1$, $m_{k-1}$ is even by Theorem 1. If we let $m_{k-1}=2m$, $x_k=2^{2m+1} n+\sum_{j=0}^{m-1} 4^j$ , where $m$ is positive integer. Since $m_{k-1}$ is odd for $x_{k-1}= 6n+5$, if we let $m_{k-1}=2m-1$, $x_k=2^{2m-1} (2n+1)+\sum_{j=0}^{m-1} 4^j$ . Each formula for the first generation odds generated from an odd includes the term $\sum_{j=0}^{m-1} 4^j$.

The $x_{k-1}$ concerns every odd excluding multiple of 3 ready for ascending operation $x_k=\left(2^{m_{k-1}} x_{k-1} -1\right)/3$, and the $x_k$ concerns every odd ready for descending operation $x_{k-1}=\left(3x_k +1\right)/2^{m_{k-1}}$. Thus, descending operation against any odds can generate the odds not multiple of 3, and vice versa. Descending operations against the odds  1,5,21, $\cdots \left(= \sum_{j=0}^{m} 4^j \right)$  reach the origin $x_0=1$, and otherwise the odds greater than or equal to 5. This proves Theorem 3.

\begin {theorem}
Let $x_k=2^{2m+1} n+\sum_{j=0}^{m-1} 4^j$  be the odds generated by $x_k=\left(2^{2m} x_{k-1} -1\right)/3$ operation against odds $x_{k-1}=6n+1$ and $x_k=2^{2m-1} (2n+1)+\sum_{j=0}^{m-1} 4^j$  be the odds generated by $x_k=\left(2^{2m-1} x_{k-1}-1\right)/3$ operation against odds $x_{k-1}=6n+5$, then $x_{k-1}$ concern every odd excluding multiple of 3 and $x_k$ every odd, where $n$ is nonnegative integer.
\end {theorem}

\subsection{Ascending operation against 1}

If we let $x_{k-1}=x_k$ in the recurrence relation $x_k=\left(2^{m} x_{k-1} -1\right)/3$, $(2^m-3) x_{k-1}=1$. Solution of the formula is $x_{k-1}=1$ and $m=2$. This is the trivial cycle. Steiner proved that there is no 1-cycle other than the trivial (1; 2).\cite{trivial cycle} Since the trivial cycle is unique, no pairs of $x_{k-1}$ and $x_k$ are the same except $x_{k-1}=x_k=1.$ 

An ascending sequence of $x_0=x_1=x_2=\cdots=1$ is generated by the recurrence relation $x_k=\left(2^2 x_{k-1}-1\right)/3$ of ascending operation against $x_0=1$. Every $x_{k-1}=1$ in each generation generates the odds $x_k=\sum_{j=0}^{(m-2)/2} 4^j$ in the next generation. Every odd $x_k$ in the sequence generated from  $x_{k-1}=1$ goes to the origin $x_0=1$ by descending operations. The odd $x_{k-1}=1$ for $k\geq2$ can be thus considered exceptionally.

\subsection{Ascending operation against odds}

Since $2^3\equiv -1 \left(mod\; 3^2 \right)$ and $2^6\equiv 1 \left(mod\; 3^2 \right)$, $2^{6n}\equiv 1 \left(mod\; 3^2 \right)$, where $n$ is positive integer. Let $x_0$ be a positive odd not multiple of 3, then $2^{6n+m_0 } x_0\equiv 2^{m_0 } x_0  \left(mod \; 2\times3^2 \right)$ by the properties of modular arithmetic regarding compatibility with scaling, where $m_0$ is nonnegative integer less than 6. 

Let $x_1=\left(2^{m_0} x_0-1\right)/3$ and $y_1=\left(2^m x_0-1\right)/3$ be odds with $m\equiv m_0  \left(mod\; 6\right)$, then $y_1\equiv x_1  \left(mod\; 6\right)$. This means that the first generation odds in ascending operation are congruent modulo 6, if numbers of doubling operations are congruent modulo 6, as illustrated in Tables 1 and 2. This proves Theorem 4. 

\begin {theorem}
Let $x_0$ be positive odds not multiple of 3, and let $x_1=\left(2^{m_0} x_0-1\right)/3$ and $y_1=\left(2^m x_0-1\right)/3$ be odds with $m\equiv m_0  \left(mod \; 6\right)$, then $y_1\equiv x_1  \left(mod\; 6\right)$, where $m_0$ is nonnegative integer less than 6. 
\end {theorem}

An ascending sequence of odds from $x_0$ to $x_n$ generated by the recurrence relation $x_k=\left(2^m x_{k-1}-1\right)/3$ of ascending operation can be expressed as $x_0 \Rightarrow x_1 \Rightarrow\cdots\Rightarrow x_k \Rightarrow \cdots \Rightarrow x_n$. If $3\mid x_n$, the $x_n$ is a terminal number of ascending sequence because  no further ascending operation is possible. The odds $x_0$ to $x_{n-1}$ are not multiples of 3.

Let an odd $x_1=\left(2^{m_0} x_0-1\right)/3$ such that $x_1 \equiv 1 \left(mod \; 3 \right)$, then $\left(2^{m_0+2} x_0-1\right)/3\equiv -1 \left(mod \; 3\right)$ and $\left(2^{m_0+4} x_0-1\right)/3\equiv 0 \left(mod \; 3\right)$. Let $y_1=\left(2^m x_0-1\right)/3$, then $y_1\equiv1 \left(mod\; 3\right)$ for $m\equiv m_0  \left(mod\; 6\right)$, $y_1\equiv -1 \left(mod\; 3\right)$ for $m\equiv m_0+2 \left(mod\; 6\right)$ and $y_1\equiv 0 \left(mod\; 3\right)$ for $ m\equiv m_0+4 \left(mod\; 6\right)$ by Theorem 4. The odds multiple of 3 are generated by every 6 doubling operations of ascending operation. Since the first generation odds are generated by every two doubling operations, ratio of the odds multiple of 3 against the odds in the first generation is 1/3. 

Infinite number of odds multiple of 3 are generated by one step of ascending operation against any positive odds not multiple of 3. And reversely, odds not multiple of 3 are generated by one step of descending operation against any positive odds multiple of 3. Every positive odd not multiple of 3 is thus generated by the descending operation against every positive odd multiple of 3, and \emph{vice versa}. Therefore, every descending sequence can start from terminal numbers. 

\subsection{Odds in the branches}

Number of odds $x_1$ in the branch of first generation originated from an odd $x_0$ is infinite. The odds $x_k$ in the branches of $k$-th generation are infinitely generated by ascending operations against each odd $x_{k-1}$ excluding multiples of 3 in the branches of previous generation. The sequence of branches can be thus generated infinitely by  consecutive ascending operations. 

If we let $M$ be number of odds excluding multiples of 3 in each branch of every generation, number of odds excluding multiples of 3 in the branches of $k$-th generation is $M^k$, because number of odds ready for ascending operation in the branches of previous generation is $M^{k-1}$. Total number of odds excluding multiples of 3 up to branches of $n$-th generation from an odd is $1+M+M^2+\cdots+M^n$. If we let $1 \ll M$, since $M^{n-1} \ll M^n$, $1+M+M^2+\cdots +M^n \approx M^n$. If the odds in in every branch are not limited, $M$ and $n$ go to infinite. 

\section {Generalization for patterns of ascending operation}
\subsection {Modular arithmetic properties of doubling operation}

Let $p$ be an odd prime number, then $\left(p\pm1\right)^{p^n} = \sum_{k=0}^{p^n} \left( \begin{array}{c} p^n \\ k \end{array} \right) \left(\pm1\right)^k p^{p^n-k}$, where $n$ is nonnegative integer. Since p-adic valuation for the binomial coefficient is $\nu_p \left( \left( \begin{array}{c} p^n \\ k \end{array} \right) \right) = n- \nu_p (k)$, $\nu_p \left( \left( \begin{array}{c} p^n \\ k \end{array} \right) p^{p^n-k} \right) = n-\nu_p (k) + p^n-k$. Since $k+\nu_p (k)< p^n$ if $0\leq k\leq p^n-1$ and $n-\nu_p (k)+p^n-k=0$ if $k=p^n$, $p^{n+1}\mid \left( \begin{array}{c} p^n \\ k \end{array} \right) p^{p^n-k}$ for $0\leq k\leq p^n-1$. Thus, $\left(p\pm1\right)^{p^n} \equiv \pm 1 \left( mod\; p^{n+1} \right)$. This proves Theorem 5. The prime number 2 satisfies also the formula $\left(p+1\right)^{p^n} \equiv 1 \left( mod\; p^{n+1} \right)$, thus $3^{2^n} \equiv 1 \left( mod\; 2^{n+1} \right)$. 

\begin {theorem}
Let $p$ be an odd prime number, then $\left(p\pm1\right)^{p^n} \equiv \pm 1 \left( mod\; p^{n+1} \right)$, where $n$ is nonnegative integer. 
\end {theorem}

If we let $p=3$, $4^{3^n}=2^{2\times 3^n}\equiv 1 \left( mod\; 3^{n+1} \right)$ and $2^{3^n} \equiv -1 \left( mod\; 3^{n+1} \right)$ by Theorem 5. Since the odds $2^{3^n}-1$ and $2^{3^n}+1$ in $2^{2\times 3^n}-1= \left( 2^{3^n}-1\right) \left( 2^{3^n}+1\right)$ are neighboring, $3^{n+1} \mid 2^{3^n}+1$ and $3 \nmid 2^{3^n}-1$. Mersenne numbers show $3^{n+1} \mid 2^{2\times 3^n}-1$ and $3 \nmid 2^{3^n}-1$, and Wagstaff numbers show $3^n \mid \left( 2^{3^n}+1 \right)/3$. This proves Theorems 6 and 7.

\begin {theorem}
Mersenne numbers show $3^{n+1} \mid 2^{2\times 3^n}-1$ and $3 \nmid 2^{3^n}-1$, where $n$ is nonnegative integer.
\end {theorem}

\begin {theorem}
Wagstaff numbers show $3^n \mid \left( 2^{3^n}+1 \right)/3$, where $n$ is nonnegative integer.
\end {theorem}

 Since $3^{n+1} \mid 2^{3^n}+1$ and $3 \nmid 2^{3^n}-1$ for $2^{2\times 3^n}-1= \left( 2^{3^n}-1\right) \left( 2^{3^n}+1\right)$, $3^n \mid \sum_{j=0}^{3^n-1} 4^j$, $3 \nmid \sum_{j=0}^{3^n-1} 2^j$, $3^n \mid \sum_{j=0}^{3^n-1} \left(-2 \right)^j$  and $\sum_{j=0}^{3^n-1} 4^j = \sum_{j=0}^{3^n-1} 2^j\times\sum_{j=0}^{3^n-1} \left(-2 \right)^j$.  Generally, since $x^{2m}-1= \left( x^m-1\right) \left( x^m+1\right)$  for positive number $x$ and odd $m$, $\sum_{j=0}^{m-1} x^{2j} = \sum_{j=0}^{m-1} x^j\times\sum_{j=0}^{m-1} \left(-x \right)^j$. 

Since $x^{2^{a}m}-1 = \left( x^m-1 \right)\prod_{k=0}^{a-1} \left( x^{2^{k}m}+1 \right)$ and $x^{2^a}-1 = \left( x-1 \right)\prod_{k=0}^{a-1} \left( x^{2^k}+1 \right)$, $\sum_{j=0}^{2^a-1} x^j = \prod_{k=0}^{a-1} \left( x^{2^k}+1 \right)$ and $\sum_{j=0}^{m-1} x^{2^{a}j} = \sum_{j=0}^{m-1} x^j \times \prod _{k=0}^{a-1} \sum_{j=0}^{m-1} \left( -x^{2^k} \right)^j$. If we let x be even, $x^{2^k}+1$ with $k\geq1$ are generalized Fermat numbers, This proves Theorem 8. 

\begin {theorem}
$\sum_{j=0}^{2^a-1} x^j = \prod_{k=0}^{a-1} \left( x^{2^k}+1 \right)$ and $\sum_{j=0}^{m-1} x^{2^{a}j} = \sum_{j=0}^{m-1} x^j \times \prod _{k=0}^{a-1} \sum_{j=0}^{m-1} \left( -x^{2^k} \right)^j$ for positive number $x$ and odd $m$.
\end {theorem}

Theorem 6 states $2^{2\times3^n}\equiv 1 \left( mod \; 3^{n+1} \right)$, where $n$ is nonnegative integer. Then, an positive odd $x_0$ not multiple of 3 shows $2^{m_0 +2\times 3^n} x_0 \equiv 2^{m_0} x_0  \left( mod \; 2\times 3^{n+1} \right)$ by the properties of modular arithmetic regarding compatibility with scaling, where $m_0$ is nonnegative integer less than $2\times3^n$. If we let $x_1 =\left(2^{m_0} x_0-1\right)/3$ and $y_1=\left(2^m x_0-1\right)/3$ be odds with relation of $m \equiv m_0 \left( mod \; 2\times 3^n \right)$, $y_1 \equiv x_1  \left( mod \; 2\times 3^n \right)$, where $x_0$ are odds not multiple of 3. Thus, the first generation odds generated by ascending operations are congruent modulo $2\times3^n$, if the numbers of doubling operations are congruent modulo $2\times3^n$. This proves Theorem 9. Theorem 9 is a generalized Theorem 4.

\begin {theorem}
Let $x_0$ be a positive odd not multiple of 3 and let $y_1=\left(2^m x_0-1\right)/3$ and $x_1 =\left(2^{m_0} x_0-1\right)/3$ with $m \equiv m_0 \left( mod \; 2\times 3^n \right)$ be odds, then $y_1 \equiv x_1  \left( mod \; 2\times 3^n \right)$, where $m_0$ is nonnegative integer less than $2\times3^n$ and $n$ is positive integer. 
\end {theorem}

\subsection{Pattern of ascending operations}

Let $x_0-\left(m_0 \right)-x_1-\left(m_1 \right)-x_2-(m_2 )-\cdots-x_{k-1}-\left(m_{k-1} \right)-x_k-\cdots-\left(m_{n-1} \right)-x_n$ be a sequence of odds generated by the recurrence relation $x_k = \left(2^{m_{k-1}} x_{k-1}-1 \right)/3$ of ascending operation. Each number in the parenthesis is number of doubling operations. The ascending sequence can be expressed as $x_0 \Rightarrow x_1 \Rightarrow \cdots \Rightarrow x_{k-1} \Rightarrow x_k \Rightarrow \cdots \Rightarrow x_n$. Pattern of a sequence is identified by the numbers of doubling operations in the sequence ignoring the sizes of odds.

Let $b_k=\sum_{h=0}^{k-1} m_h$  be the total number of doubling operations in a sequence from $x_0$ to $x_k$. If we let $y_0=x_0+2\times 3^n j$, odds in the ascending sequence $y_1= x_1+ 2^{m_0+1} 3^{n-1} j$, $y_2= x_2+ 2^{m_0+ m_1 +1} 3^{n-2} j$, $y_k=x_k+ 2^{b_k+1} 3^{n-k} j$ and $y_n=x_n+2^{b_n+1} j$ are generated by the recurrence relation $y_k = \left(2^{m_{k-1}} y_{k-1}-1 \right)/3$ of ascending operation, where $j$ is nonnegative integer. The ascending sequences $y_0 \Rightarrow  y_1 \Rightarrow y_2 \Rightarrow \cdots \Rightarrow y_k \Rightarrow \cdots \Rightarrow y_n$ with different $j$ are generated by the same numbers of doubling operations as the sequence from $x_0$ to $x_n$ in the same pattern. The sequences in the same pattern show $y_0\equiv x_0  \left( mod\; 2\times 3^n \right)$, $y_1\equiv x_1  \left( mod\; 2^{m_0 +1} 3^{n-1} \right)$, $y_k\equiv x_k  \left( mod\; 2^{b_{k+1}} 3^{n-k} \right)$, and $y_n\equiv x_n  \left( mod\; 2^{b_n +1} \right)$. Number of sequences in the same pattern is infinite. This proves Lemma 10. 

Each descending sequence started from each odd is unique. Thus $y_0 = x_0 + 2\times 3^n j$ is the solution for the recurrence relation $x_{k-1}=\left(3x_k+1\right)/2^{m_{k-1}}$ of descending operation against $y_n = x_n + 2^{b_n+1}j$. The descending sequences from $y_n$ to $y_0$ with different $j$ belong to the same pattern. 

\begin {lemma}
Let $x_0 \Rightarrow x_1 \Rightarrow x_2 \cdots \Rightarrow x_k \Rightarrow \cdots \Rightarrow x_n$ be a sequence generated by the recurrence relation $x_k = \left(2^{m_{k-1}} x_{k-1}-1 \right)/3$ of ascending operation with $b_k=\sum_{h=0}^{k-1} m_h$, if and only if $y_0\equiv x_0  \left( mod\; 2\times 3^n \right)$, then the sequences of odds $y_0 \Rightarrow  y_1 \Rightarrow y_2 \Rightarrow \cdots \Rightarrow y_k \Rightarrow \cdots \Rightarrow y_n$ with $y_k\equiv x_k  \left( mod\; 2^{b_{k+1}} 3^{n-k} \right)$ are infinitely generated by the recurrence relation $y_k = \left(2^{m_{k-1}} y_{k-1}-1 \right)/3$ with the same numbers of doubling operations $m_k$.
\end {lemma}

\subsection{Primitive sequence}

Theorem 9 states that let $x_0$ be positive odds not multiple of 3, then the first generation odds $y_1 = \left(2^{m} x_0-1 \right)/3$ and $x_1 = \left(2^{m_0} x_0-1 \right)/3$ with $m\equiv m_0  \left( mod\; 2\times 3^n \right)$ show $y_1\equiv x_1  \left( mod\; 2\times 3^n \right)$, where $m_0$ is nonnegative integer less than $2\times3^n$ and $n$ is positive integer. Then, $n$ terms of ascending sequences started from each $y_1$ belong to the same pattern by Lemma 10 

Let a group such that $G_n = \left\{x_1 \mid 1\leq x_1< 2 \times 3^n, \left( 2^m x_0 -1 \right) /3 \equiv x_1 \left(mod\; 2\times 3^n \right), 2^m x_0 \equiv 1 \left(mod\; 3 \right), m\equiv m_0 \left(mod\; 2\times 3^n \right), 0\leq m_0< 2 \times 3^n, m_0 \in \mathbb{N},n\in \mathbb{N} \right\}$. Order of the group is $\vert G_n \vert =3^n$, since $2^m x_0 \equiv 1 \left( mod\; 3\right)$ by every two doubling operations. Since $\left(2^m x_0-1\right)/3\equiv \left(2^{m_0} x_0-1\right)/3 \left(mod\; 2\times 3^n \right)$ by Theorem 9,  $G_n=\left\{ x_1 \mid 1\leq x_1< 2 \times 3^n, \left(2^{m_0} x_0 -1\right)/3 \equiv x_1 \left(mod\; 2\times 3^n \right), 2^{m_0} x_0\equiv 1 \left(mod\; 3 \right), 0\leq m_0< 2\times 3^n, m_0 \in \mathbb{N}, n\in \mathbb{N} \right\}$. Since $G_n$ includes every positive odd less than $2\times 3^n$, $G_n= \left\{ x \mid x=2m-1, 1\leq m\leq 3^n, m\in \mathbb{N}, n\in \mathbb{N} \right\}$.

For example, $G_2=\left\{x_1 \mid 1\leq x_1< 2 \times 3^n, \left(2^{m+2\times 3^2 j}-1 \right)/3\equiv x_1 \; \left( mod\; 2\times 3^2 \right), 2^{m+2\times 3^2 j} \equiv 0\; \left( mod\; 3\right), 0\leq m< 2\times 3^2, m_0\in \mathbb{N}, n\in \mathbb{N} \right\}$. Since $\left( 2^{m+2\times 3^2 j}-1 \right) /3 \equiv \; \left( 2^m-1\right)/3 \; \left( mod \; 2\times 3^2 \right)$ by Theorem 9, $G_2 = \left\{ x_1 \mid 1\leq x_1< 2 \times 3^n, \left(2^m -1\right)/3 \equiv x_1 \left(mod\; 2\times 3^n \right), 2^m \equiv 1 \left (mod \; 3 \right), 0\leq m< 2\times 3^n, m \in \mathbb{N} \right\}= \left\{ 1,5,3,13,17,15,7,11,9 \right\}$. The elements of the group are remainders of $y_1\equiv x_1  \left( mod\; 2\times 3^n \right)$. Thus, $G_2= \left\{ x \mid x=2m-1, 1 \leq m\leq 3^n, m \in \mathbb{N} \right\}$. 

The $n$ terms of ascending sequences started from each $y_1$ with relation of $y_1\equiv x_1  \left( mod\; 2\times 3^n \right)$ show the same pattern as the sequence started from $x_1$. If $ x_1 \in G_n$, the ascending sequence started from $x_1$ is named as a primitive sequence, which  is generated by the smallest odds among the sequences in the same pattern. Thus, each set of ascending sequences in the same pattern includes each primitive sequence that starts from an odd between 1 and $2\times3^n$. The sequences in the same patterns are infinitely regenerated throughout the consecutive generations like a fractal.

\section{Specific consecutive descending operations}
\subsection{Descending operations with one halving operation}

An odd $x_{k-1}$ generated by the recurrence relation $x_{k-1}=\left( 3x_k+1\right) /2^{m_{k-1}}$ of descending operation is smaller than the odd $x_k$ if $m_{k-1}\geq 2$. If average numbers of halving operations in a sequence is $\overline{m_{k-1}}\approx \log_2 3 \approx 1.59$, the sequence will keep the similar values. If $m_{k-1}=1$, $x_{k-1}>x_k$. 

Let $x_n-[1]-x_{n-1}-[1]-x_{n-2}-[1]-\cdots x_k-\cdots -[1]-x_0$ be a sequence  generated by the recurrence relation $x_{k-1}=\left( 3x_k+1\right) /2$ of descending operation with one halving operation. Numbers of halving operations are listed in the brackets. Total number of halving operations is $b_n=n$.  Sizes of the odds in the descending sequence increase continuously.

Odd 3 is generated by the ascending operation $3=\left(2\times5-1\right)/3$ against 5. Then, ascending operations with one doubling operation against the odds $x_0=-1+2\times3j$ generate the odds $x_1 =-1+4j$ by Lemma 10, where j is posiive integer. Two terms of sequences with $m_0=m_1=1$ are $x_2=-1+2^3  j$, $x_1=-1+2^2 3j$ and $x_0=-1+2\times 3^2 j$. Three terms of sequences with $m_0=m_1=m_2=1$ are $x_3=-1+2^4 j$, $x_2=-1+2^3 3 j$, $x_1=-1+2^2 3^2 j$ and $x_0=-1+2\times3^3 j$. Thus, the odds $x_k=-1+2^{k+1} 3^{n-k} j$ belong to the $n$ terms of descending sequences from $x_n=-1+2^{n+1} j$ to $x_0=-1+2\times 3^n j$. The descending sequences of different $j$ show the same pattern.

Since $x_0/x_n\approx (3/2)^n$, sizes of the odds increase exponentially to the number of descending operations. The ratio $x_0/x_n$ becomes more close to $(3/2)^n$ by larger $x_n$. 

Let $x_{n+1}-[1]-x_n-[1]-x_{n-1}-[1]-\cdots -x_2-[1]-x_1-[m_0 ]-1$ be a descending sequence, where $x_1$ is a first generation odd originated from $x_0=1$. The odds $x_{n+1}$ to $x_1$ increase geometrically with common ratio of $x_{k-1}/x_k\approx 3/2$. 

Theorem 7 states $3^{n+1}\mid 2^{3^n}+1$, where $n$ is nonnegative integer. Since $3^{n+1} \mid 2^{3^n (2h-1)}+1$, $\left(2^{3^n (2h-1)+1}-1\right)/3 \equiv -1 \left( mod\; 2\times 3^n \right)$, where $h$ is positive integer. If we let $x_1=\left(2^{3^n (2h-1)+1}-1\right)/3$, since $x_1=\sum _{j=0}^{m_0 -1} 4^j$  with $m_0 =\left(3^n \left(2h-1\right)+1\right)/2$, the $x_1$ belong to the first generation odds generated from the origin $x_0 =1$ as shown in Theorem 3. Since $x_1 \equiv -1 \left( mod\; 2\times 3^n \right)$, $n$ terms of ascending sequences $x_1-(1)-x_2-(1)-\cdots-x_{n}-(1)-x_{n+1}$ in the same pattern are generated by Lemma 10. If we let $x_1=-1+2\times 3^n j$, the sequence of odds from $x_1$ to $x_{n+1}=-1+2^{n+1} j$ with $x_k=-1+2^k 3^{n-k+1} j$ are generated by the recurrence relation $x_k=\left( 2 x_{k-1} -1\right)/3$, where $j$ is positive integer.  

Number of halving operations to reach the origin $x_0=1$ is $m_0=3^n \left(2h-1\right)+1$. Infinitely large $x_1$ can be generated by consecutive descending operations with one halving operation. One step of descending operation against the $x_1$ reaches the origin $x_0 =1$ irrespective of size of $x_1$. If we let $h=1$, $m_0 = 3^n+1$, such as 4, 10, 28, $\cdots$.

Examples of descending sequences of odds with $h=1$ are $3-[1]-5-[4]-1$ for $n=1$, Then, $3\equiv -1\; \left( mod\; 2^2 \right)$ and $5\equiv -1\; \left( mod\; 2\times3 \right)$. If $n=2$, $151-[1]-227-[1]-341-[10]-1$, Then, $151\equiv -1\; \left( mod\; 2^3 \right)$, $227\equiv -1\; \left( mod\; 2^2 3 \right)$ and $341\equiv -1\; \left( mod\; 2\times 3^2 \right)$. If $n=3$, $26512143-[1]-39768215-[1]-59652323-[1]-89478485-[28]-1$. Then,  $26512143\equiv -1\; \left( mod\; 2^4 \right)$, $39768215\equiv -1\; \left( mod\; 2^3 3 \right)$, $59652323 \equiv -1\; \left( mod\; 2^2 3^2 \right)$ and $89478485\equiv -1\; \left( mod\; 2\times 3^3 \right)$. 

\subsection{Descending operations with two halving operations}

Let $x_n-[2]-x_{n-1}-[2]-x_{n-2}-[2]-\cdots-x_1-[2]-x_0$ be a sequence generated by the recurrence relation $x_{k-1}=\left( 3x_k+1\right)/2^2$ of descending operation with two halving operations. Total number of halving operations is $b_n=2n$. Sizes of the odds in the descending sequence decrease continuously.

Odd 1 is generated by the ascending operation $1=\left(2^2\times1-1\right)/3$ against 1. Then, ascending operations against the odds $x_0 =1+2\times3j$ generate the odds $x_1=1+2^3 j$ by Lemma 10, where j is nonnegative integer. Two terms of sequences with $m_0=m_1=2$ are $x_2=1+2^5 j$, $x_1=1+2^3 3 j$ and $x_0=1+2\times 3^2 j$. Three terms of sequences with  $m_0=m_1=m_2=2$ are $x_3=1+2^7 j$, $x_2=1+2^5 3 j$, $x_1=1+2^3 3^2 j$ and $x_0=1+2\times 3^3 j$. The odds $x_k=1+2^{2k+1}3^{n-k} j$ belong to the $n$ terms of descending sequences from $x_n=1+2^{2n+1} j$ to $x_0=1+2\times 3^n j$. The descending sequences of different $j$ show the same pattern.

Since $x_0/x_n\approx\left(3/4\right)^n$, sizes of the odds decrease exponentially to the number of descending operations with two halving operations. If we let $x_n=1$, $1-[2]-1-[2]-\cdots-1$, which is the $n$ terms of descending sequence concerning the trivial cycle. 

Let $x_{n+1}-[2]-x_n-[2]-x_{n-1}-[2]-\cdots-x_2-[2]-x_1-[m_0 ]-1$ be a descending sequence, where $x_1$ is a first generation odd originated from $x_0=1$. The odds $x_{n+1}$ to $x_1$ decrease geometrically with common ratio of $x_{k-1}/x_k\approx3/4$. 

Theorem 6 states $3^{n+1} \mid 2^{2\times3^n }-1$, where $n$ is nonnegative integer. Since $2^{2\times 3^n h}\equiv 1 \left( mod\; 3^{n+1} \right)$ ), $\left( 2^{2\left(3^n h+1\right)} -1\right)/3\equiv 1 \left( mod\; 2\times 3^n \right)$, where $h$ is positive integer. If we let $x_1=\left( 2^{2\left(3^n h+1\right)} -1\right)/3$, since $x_1=\sum _{j=0}^{3^n h} 4^j$, $x_1$ belong to the first generation odds generated from the origin $x_0 =1$ as shown in Theorem 3. Since $x_1\equiv 1 \left( mod\; 2\times 3^n \right)$, $n$ terms of sequences $x_1-(2)-x_2-(2)-\cdots-(2)-x_{n+1}$ in the same pattern are generated by Lemma 10. If we let $x_1=1+2\times3^n j$, sequences of odds from $x_1$ to $x_{n+1}=1+2^{2n+1}j$ with $x_k=1+2^{2k-1} 3^{n-k+1} j$ are generated by the recurrence relation $x_k=\left( 2^2 x_{k-1} -1\right)/3$, where $j$ is positive integer. Number of halving operations against the odd $x_1$ to reach the origin $x_0=1$ is $m_0=2\left(3^n h+1\right)$. If we let $h=1$, $m_0 = 2 \left( 3^n +1 \right)$, such as 2, 8, 20,56, $\cdots$.

Examples of descending sequences for $h=1$ are: $1-[2]-1$ for $n=0$. $113-[2]-85-[8]-1$ if $n=1$. Then, $113\equiv 1\; \left( mod\; 2^3 \right)$ and $85\equiv 1\; \left( mod\; 2\times 3 \right)$. If $n=2$, $621377-[2]-466033-[2]-349525-[20]-1$. Then, $621377\equiv 1\; \left( mod\; 2^5 \right)$, $466033\equiv 1\; \left( mod\; 2^3 3 \right)$ and $349525\equiv 1\; \left( mod\; 2\times 3^2 \right)$. 

\section{Sequence of odds repeating cyclically}
\subsection{Sequence including cyclic odds}

Let $\cdots-x_a-\left[m_{a-1} \right]-x_{a-1}-[m_{a-2} ]-x_{a-2}-\cdots-x_{b+1}-[m_b ]-x_b-\left[m_{b-1} \right]-\cdots$ be a descending sequence of odds that are not multiples of 3. Formula for the odds concerning $a-b$ times of descending operations is $x_b=\left( 3 \cdots \left( 3 \left( 3 \left( 3x_a+1 \right)/2^{m_{a-1}}+1\right)/2^{m_{a-2}}+1\right)/2^{m_{a-3}}\cdots +1\right)/2^{m_b}$.

Lemma 2 states no two odds $x_k$ in a generation are the same, if the odds $x_{k-1}$ in the previous generation are different. If we let $x_a = x_b$ in a sequence, since the descending sequence started from an odd is unique, $x_{a-1} = x_{b-1}$, $x_{a-2} = x_{b-2}$, $x_{a-3} = x_{b-3}$, $\cdots$ and $m_a = m_b$, $m_{a-1} =m_{b-1}$, $m_{a-2} =m_{b-2}$, $\cdots$. The sequence is consisted with the cyclic odds from $x_a$ to $x_{b+1}$. 

The cyclic descending operations continue sequentially until the descending sequence reaches an origin of $x_0$. All of the odds $x_a$ to $x_{b+1}$ in the sequence have the chance to be the origin of the sequence. If the origin of the descending sequence is unique, the odds in the cyclic sequence will show $x_a=x_{a-1}=x_{a-2}=\cdots=x_{b+1} =\cdots=x_0$. Since no pairs of $x_{k-1}$ and $x_k$ are the same except $x_{k-1}=x_k=1$, it is expected that the odds in the cyclic sequence are equally 1.

\subsection{Cycle with one odd}

The cyclic sequence with one odd is $\cdots -x_a-\left[ m_{a-1} \right]-x_{a-1} \left( =x_a \right)-\left[ m_{a-1} \right]-\cdots-x_0 \left( =x_a \right)$. Recurrence relation for one descending operation with $x_a =x_{a-1}$ repeating in a cycle is $x_a=\left(3x_a+1\right)/2^{m_{a-1}}$. Solution is $x_a=1$ and $m_{a-1}=2$, and the repeating sequence is named as trivial cycle, $1-[2]-1$.\cite{trivial cycle}  

\subsection{Cycle with two odds}

The cyclic sequence with two odds is $\cdots-[m_a ]-x_a-\left[ m_{a-1} \right]-x_{a-1}-\left[ m_a \right]-x_{a-2} (=x_a )-\left[ m_{a-1} \right]-x_{a-3} \left( =x_{a-1} \right)-\left[ m_a \right]-\cdots$. Formula for the odds concerning two times of descending operations repeating in a cycle for $x_a=x_{a-2}$ is $\left(2^{m_a+m_{a-1}}-3^2 \right) x_a =2^{m_{a-1}}+3$, and $\left(2^{m_a+m_{a-1}}-3^2 \right) x_{a-1} = 2^{m_a}+3$ for $x_{a-1}=x_{a-3}$. Since $2^{m_a+m_{a-1}}>3^2$, $m_a+m_{a-1}\geq 4$. 

Since $x_{a-1}=\left(3x_a+1\right)/2^{m_{a-1}}$ and $x_a=\left(3x_{a-1}+1\right)/2^{m_a}$ are the recurrence relations of the cycle, $\left( 3+ 1/x_a \right)\left(3+1/x_{a-1} \right)=2^{m_a+m_{a-1}}$. The solution is $x_a=x_{a-1}=1$ and $m_a=m_{a-1}=2$. Thus, there is no cycle composed of two odds, but $1-[2]-1-[2]-1$. 

\subsection{Cycle with three odds}

The cyclic sequence with three odds is $\cdots-[m_a ]-x_a-\left[ m_{a-1} \right]-x_{a-1}-\left[ m_{a-2} \right]-x_{a-2}-\left[ m_a \right]-x_{a-3} (=x_a ) -\left[ m_{a-1} \right]-x_{a-4}\left( =x_{a-1} \right)-\left[ m_{a-2} \right]-x_{a-5}\left( =x_{a-2} \right)-\left[ m_a \right]-\cdots$. Formula for the odds concerning three times of descending operations repeating in a cycle for $x_a=x_{a-3}$ is $\left(2^{m_a+m_{a-1}+m_{a-2}}-3^3 \right) x_a =2^{m_{a-1}+m_{a-2}}+ 2^{m_{a-1}}3+3^2$. Since  left side of the formula is positive, $m_a+m_{a-1}+m_{a-2}\geq5$.

Recurrence relations for the descending operations of the cycle are $x_{a-1}=\left(3x_a+1\right)/2^{m_{a-1}}$, $x_{a-2}=\left(3x_{a-1}+1\right)/2^{m_{a-2}}$ and $x_a =\left(3x_{a-2}+1\right)/2^{m_a}$. Multiplication of the three formulae results $\left( 3+ 1/x_a \right) \left(3+1/x_{a-1} \right) \left(3+1/x_{a-2} \right) =2^{m_a+m_{a-1}+m_{a-2}}$. The solution is $x_a=x_{a-1}=x_{a-2}=1$ and $m_a=m_{a-1}=m_{a-2}=2$. There is no cycle composed of three odds, but $1-[2]-1-[2]-1-[2]-1$.

\subsection{Cycle with n odds}

The cyclic sequence with $a-b\; \left(=n\right)$ odds is $\cdots-[m_b ]-x_a \left( =x_b \right)-\left[ m_{a-1} \right]-x_{a-1}-\left[ m_{a-2} \right]-x_{a-2}-\cdots-x_{b+1}-\left[ m_b \right]-x_b \left( =x_a \right)-\left[ m_{a-1} \right]-\cdots$. Formula for the cyclic sequence with $x_a=x_b$ is $\left(2^{m_{a-1}+m_{a-2}+\cdots+m_b}-3^{a-b} \right) x_a =2^{m_{a-1}+m_{a-2}+\cdots+m_{b+1}}+2^{m_{a-1}+m_{a-2}+\cdots+m_{b+2}}3+\cdots+ 2^{m_{a-1}}3^{a-b-2}+3^{a-b-1}$. 

Recurrence relations for the descending operations of the cycle are $x_{a-1}=\left(3x_a+1\right)/2^{m_{a-1}}$, $x_{a-2}=\left(3x_{a-1}+1\right)/2^{m_{a-2}}$, $x_{a-3}=\left(3x_{a-2}+1\right)/2^{m_{a-3}}$,$\cdots$ , $x_a=\left(3x_{b+1}+1\right)/2^{m_b}$. Multiplication of the the formulae results $\left( 3+ 1/x_a \right) \left(3+1/x_{a-1} \right)\cdots \left(3+1/x_{b+1} \right) = 2^{m_{a-1}+m_{a-2}+\cdots+m_b}$. The solution is $x_a=x_{a-1}=\cdots=x_{b+1}=1$ and $m_{a-1}=m_{a-2}=\cdots= m_b=2$. There is no cycle composed of {n} odds, but $1-[2]-1-[2]-\cdots-1$ with $n$ terms..

\subsection{Cycle with infinite odds}

Formula for $x_a=1$ concerning $n$ times of operation with $m_k=2$ is $2^2-2^2 \left(3/2^2 \right)^n =1+3/2^2+\left(3/2^2 \right)^2+\cdots +\left(3/2^2 \right)^{n-1}$. If n goes to infinite, $4=1+3/2^2+\left(3/2^2 \right)^2+\left(3/2^2 \right)^3+\cdots$. 

It is concluded that no odd but 1 is the number in the cyclic sequence. Thus, no two odds except 1 are the same in descending sequences. This proves Lemma 11. 

\begin {lemma}
Let a sequence of positive odds $x_n\rightarrow x_{n-1}\rightarrow \cdots \rightarrow x_k \rightarrow x_{k-1} \rightarrow \cdots \rightarrow x_1\rightarrow x_0$ generated by the recurrence relation $x_{k-1}=\left(3x_k+1\right)/2^{m_{k-1}}$ of descending operation with $m_{k-1}=\nu_2 \left(3x_k+1\right)$, then no two odds except 1 in the sequence are the same.
\end {lemma} 

\section{Estimation of odds in descending operations}

Approximation for the ratio of neighboring odds $x_{k-1}/x_k\approx3/2^{m_{k-1}}$ is derived from the recurrence relation $x_{k-1}=\left(3x_k +1\right)/2^{m_{k-1}}$. Since the chain relation shows $\left(x_{n-1}/x_n \right)\left(x_{n-2}/x_{n-1} \right)\cdots\left(x_0/x_1 \right)\approx 3^n/2^{b_n}$ with $b_n=\sum_{j=0}^{n-1} m_j$, $x_0\approx 3^n x_n/2^{b_n}$. If we let $x_0=1$ in a descending sequence of odds, total number of halving operations for n times of descending operations against an odd $x_n$ to reach the origin of 1 is $b_n\approx \log_2 {3^n x_n}$. 

Lemma 10 states that the descending sequences from $y_n=x_n+2^{b_n+1} j$ to $y_0=x_0 +2\times 3^n j$ generated by the recurrence relation $y_{k-1}=\left(3y_k +1\right)/2^{m_{k-1}}$ show the same pattern as the sequence from $x_n$ to $x_0$, where $j$ is nonnegative integer.The ratio $y_n /y_0 \approx2^{b_n}/ 3^n$ of the sequences in the same pattern is approximately constant, since the sequences from $y_n$ to $y_0$ have the same $b_n$. The ratio of $y_n /y_0$ by larger odds becomes more close to $2^{b_n}/ 3^n$. This proves Lemma 12.

\begin {lemma}
Let $b_n=\sum_{k=0}^{n-1} m_k$  be the total number of doubling operations in the sequence of odds $x_n\rightarrow x_{n-1}\rightarrow \cdots \rightarrow x_k \rightarrow \cdots \rightarrow x_0$ generated by the recurrence relation of descending operation $x_{k-1}=\left(3x_k +1\right)/2^{m_{k-1}}$ with $m_{k-1}=\nu_2\left (3x_k+1\right)$,if $y_n\equiv x_n  \left( mod\; 2^{b_n+1} \right)$, then each sequence of odds $y_n\rightarrow y_{n-1}\rightarrow \cdots \rightarrow y_k \rightarrow \cdots \rightarrow y_0$ generated by $y_{k-1}=\left(3y_k +1\right)/2^{m_{k-1}}$ shows $y_0\approx 3^n y_n/2^{b_n}$.
\end {lemma}

An example of descending sequence of odds with one halving operation reaching $x_0=1$ is $26512143-[1]-39768215-[1]-59652323-[1]-89478485-[28]-1$ with $b_n=31$  for $n=4$. The expected value is $x_4\approx 2^{31}/3^4= 26512143.8$ for $x_0 =1$. 

General formula for the odds in descending sequence with one halving operation is $y_k=-1+2^{k+1} 3^{n-k} j$. If $n=3$, $y_3 /y_0\approx3^3 /2^3$ by Lemma 12. The primitive descending sequence of $j=1$ is $15-[1]-23-[1]-35-[1]-53$. The sequences from $x_4$ to $x_1$ and from $y_3$ to $y_0$ belong to the same pattern. Since $x_4$ is larger than $y_3$ of $j=1$, the ratio of $x_4/x_1$ is more close to $3^3 /2^3$ than $y_3/y_0$. 

\section{Origin of sequences}

The first generation odds $x_k$ are infinitely generated by the recurrence relation $x_k=\left(2^{m_{k-1}} x_{k-1}-1\right)/3$ of ascending operation with $m_{k-1}$ times of doubling operations against a positive odd $x_{k-1}$ not multiple of 3. The first generation odds show sequentially $x_k\equiv h\; (mod\; 3)$ with $h=1, -1, 0$ by every two doubling operations $m_{k-1}$. Consecutive ascending operations can be performed infinitely against the odds not multiples of 3 in the previous generation. Ascending operations are terminated by the odds multiple of 3.

Each descending sequence is unique, and can start from terminal numbers. A descending operation against any odds can generate odds not multiple of 3, and vice versa. Descending operations against the odds of 1,5,21,85,$\cdots\;\left(=\sum_{j=0}^{m} 4^j \right)$ reach 1, and otherwise reach the odds greater than or equal to 5. Since no odds in each descending sequence except 1 are the same, every descending sequence started from terminal numbers continues until to reach the origin of 1. This proves the Collatz conjecture. 

\section*{Funding}

This research did not receive any specific grant from funding agencies in the public, commercial, or not-for-profit sectors.


\bigskip

\begin{thebibliography}{9}

\bibitem{CollatzBook}
O'Connor, J.J.; Robertson, E.F. (2006). "Lothar Collatz". St Andrews University School of Mathematics and Statistics, Scotland. 

\bibitem{trivial cycle}
Steiner, R. P. (1977). "A theorem on the syracuse problem". Proceedings of the 7th Manitoba Conference on Numerical Mathematics. pp. 553–9. MR 0535032. 

\end{thebibliography}
\end{document}